\documentclass[11pt]{article}

\usepackage{amscd}
\def \text{\mbox}

\title{BGP-reflection functors and cluster combinatorics
\thanks{Supported by the NSF of China (Grants 10471071)
and by in part the Doctoral Program Foundation of Institute of Higher
Education(2003)}}

\author{Bin
Zhu\thanks{E-mail: bzhu@math.tsinghua.edu.cn}}

\date{Department of Mathematical Sciences, \\ Tsinghua University,
100084 Beijing, P. R. China}

\begin{document}

\maketitle

\def\Longrightarrow{{\longrightarrow}}
\def\P{{\cal P}}
\def\A{{\cal A}}
\def\T{\cal T}
\def\m{\textbf{ M}}
\def\t{{\tau }}
\def\b{\textbf{d}}
\def\K{{\cal K}}

\def\G{{\Gamma}}
\def\e{\mbox{exp}}

\def\righta{\rightarrow}

\def\s{\stackrel}

\def\ncong{\not\cong}

\def\mathbb{\NN}

\def\Hom{\mbox{Hom}}
\def\Ext{\mbox{Ext}}
\def\ind{\mbox{ind}}
\def\coprod{\amalg }

\begin{center}

\begin{minipage}{12cm}{\footnotesize\textbf{Abstract.}
 We define Bernstein-Gelfand-Ponomarev reflection functors in
 the cluster categories of hereditary algebras. They
  are triangle equivalences which provide a natural quiver
realization of the "truncated simple reflections" on the set of
almost positive roots $\Phi_{\ge -1}$ associated to a finite
dimensional semisimple Lie algebra. Combining with the tilting
theory in cluster categories developed in [4], we give a
unified interpretation via quiver representations for the
generalized associahedra associated to the root systems of all
Dynkin types (a simply-laced or non-simply-laced). This confirms
the conjecture 9.1 in [4] in all Dynkin types.

\medskip

  \textbf{Keywords.} BGP-reflection functor; truncated simple reflection; cluster; cluster category; compatibility degree.

\medskip

\textbf{Mathematics Subject Classification.} 16G20, 16G70, 52B11,
17B20.}
\end{minipage}
\end{center}
\medskip

\begin{center}

\textbf{1. Introduction}\end{center}

\medskip

As a model for the combinatorics of a Fomin-Zelevinsky's cluster
algebra [9, 10], the cluster category $\mathcal{C}(H)$
associated with a hereditary algebra $H$ over a field was
introduced in [4], see also [5]. It is the orbit category of the (bounded)
derived category of $H$ factored by the automorphism $G=[1]\tau^{-1}$,
where $[1]$ is the shift functor and $\tau$ the Auslander-Reiten
translation in the derived category of $H$. This orbit category is
 a triangulated category [14]. When $H$ is the path
  algebra of a quiver of Dynkin type (simply-laced case),
  it is proved in [4] that there is a
one-to-one correspondence between the set of indecomposable
objects in $\mathcal{C}(H)$ and the set of cluster variables of
the corresponding cluster algebras. This correspondence is given explicitly when the orientation of quiver is alternating,
 and under this correspondence,
tilting objects correspond to clusters. This was motivated by a
previous quiver-theoretic interpretation (using "decorated" quiver
representations) for
 generalized
 associahedra of simply-laced Dynkin type in the sense
 of Fomin-Zelevinsky [11][6], which was given in [15].

  In
the combinatorics of cluster algebras, the group of
piecewise-linear transformations of the root lattice generated by
"truncated simple reflections" $\sigma _i $ for $i\in I$ (the
index set of simple roots) plays an important role as Weyl group
in classical theory of semi-simple Lie algebra. A similar but
stronger tool in representation theory of quivers and hereditary
algebras is the so-called Bernstein-Gelfand-Ponomarev reflection
functors [2] or APR tilting functors [1]. Since, by [4],
the cluster categories provide a successful model to realize the
clusters and associahedra. It is natural to ask whether the BGP
reflection functors can be defined in the cluster categories.
These functors defined in the cluster categories should lead the
"truncated simple reflections" on the set of almost positive roots
and should be applicable to the clusters and associahedra. One of
the motivations of this work comes from [15], where the authors
gave a realization of the "truncated simple reflections" in the
category of "decorated" quiver representations. Unfortunately
their functors are not equivalences.

In this paper, we verify that it is indeed possible to define the
BGP-reflection functors in the cluster categories of hereditary
algebras (in fact they can be defined in a more general case including the case of root categories (compare [21]) ).
   The advantage of our functors (compare with [15]) is that
   the BGP-reflection
functors in cluster categories are triangle equivalences. By
applying these equivalences defined in the cluster
categories to the set of almost positive roots, we obtain a
realization of the
 "truncated simple reflections" [11]. This enables us to give in a unified
 way a quiver interpretation for generalized associahedra there.
 By using this realization, the main ingredients of constructions in Section 3
  in [11] follows without much effort from tilting theory developed in [4].
  This generalizes the main results on quiver
 interpretation for generalized associahedra of simply-laced case in [15]
  and confirms
 the conjecture 9.1. [4] in all Dynkin types.

 \medskip

 \begin{center}

\textbf{2. BGP reflection functors in orbit triangulated
 categories.}
\end{center}

\medskip

It is well-known that the orbit category $D^{b}(H)/G$ of the
derived category of a finite dimensional hereditary algebra $H$ is
 a triangulated category in which the images of triangles in
 $D^{b}(H)$ under the natural
 projection are still triangles when $G$ is an automorphism satisfying some
  specific conditions (the conditions $(g1)$, $(g2)$ below) [14].
  When $G=[1]\tau ^{-1}$, the
orbit category $D^{b}(H)/G$ is called the cluster category of $H$.
 We recall some basics on orbit triangulated categories
 from [14] and basics on the cluster categories from [4, 3]. We refer [18] for the basic reference for representation theory of algebras.

Let $\cal{H}$ be a hereditary category with Serre duality and with
finite dimensional Hom-spaces and Ext-spaces over a field $K$.
Denote by $\mathcal{D} =
    D^{b}(\cal{H})$ the bounded derived category of $\cal{H}$ with shift
functor $[1]$. For any category $\cal{E}$, we will denote by
$\ind\cal{E}$ the subcategory of isomorphism classes of
indecomposable objects in $\cal{E}$; depending on the context we
shall also use the same notation to denote the set of isomorphism
classes of indecomposable objects in $\cal{E}$. For any $T$ in $\cal{H}$, we denote the subcategory of $\cal{H}$ consisting
of direct summands of direct sums of finite many copies of $T$ by add$T$. Note that add$H$ denotes the category of projective $H-$modules.

Let $G \colon \cal{D} \to \cal{D}$ be a standard equivalence, i.e. $G$ is isomorphic to the derived tensor product
  $$ -\otimes_AX: D^b(A)\rightarrow D^b(A)$$
  for some complex $X$ of $A-A-$bimodules. We also assume that $G$ satisfies the following
properties:
\begin{itemize}
\item[(g1)]{For each $U$ in $\ind \mathcal{H}$, only a finite number
of objects $G^n U$, where $n \in \mathbf{Z}$, lie in $\ind
\mathcal{H}$.}
\item[(g2)]{There is some $N \in \mathbf{N}$ such that
$\{U[n] \mid U \in \ind \mathcal{H}, n \in [-N,N] \}$ contains a
system of representatives of the orbits of $G$ on $\ind \cal{D}$.}
\end{itemize}

We denote by $\mathcal{D}/ G$ the corresponding factor category.
The objects are by definition the $G$-orbits of objects in
$\cal{D}$, and the morphisms are given by
$$\Hom_{\mathcal{D}/G}(\widetilde{X},\widetilde{Y}) =
\oplus_{i \in \mathbf{Z}}
 \Hom_{\mathcal{D}}(G^iX,Y).$$
Here $X$ and $Y$ are objects in $\cal{D}$, and $\widetilde{X}$ and
$\widetilde{Y}$ are the corresponding objects in $\mathcal{D}/G$
 (although we shall sometimes write such objects simply as $X$ and
$Y$). The orbit category $\mathcal{D}/G$ is a Krull-Schmidt
category [4] and also a triangulated category [14]. The natural functor
$\pi \colon \mathcal{D} \to \mathcal{D}/G$ is a covering functor
of triangulated categories in the sense that $\pi$ is a covering
functor and a triangle functor [20]. The shift in $\mathcal{D}/G$
is induced by the shift in $\cal{D}$, and is also denoted by
$[1]$. In both cases we write as usual $\Hom(U,V[1]) =
\Ext^1(U,V)$. We then have
$$\Ext^1_{\mathcal{D}/G}(\widetilde{X},\widetilde{Y}) =
\oplus_{i \in \mathbf{Z}} \Ext^1_{\mathcal{D}}(G^i X,Y),$$ where
$X,Y$ are objects in $\cal{D}$ and $\widetilde{X},\widetilde{Y}$
are the corresponding objects in $\mathcal{D}/G$.

 We shall mainly concern the special choice of functor $G =[1] \tau^{-1}$,
where $\tau$ is the Auslander-Reiten translation in $\cal{D}$ and
$\mathcal{H}=H-$mod is the category of finite dimensional left modules over a finite dimensional
hereditary algebras $H$. In this case the factor category
$\mathcal{D} / G$ is called the cluster category of $H$, which is
denoted by $\mathcal{C}(H)$. It is not difficult to see that
ind$\mathcal{C}(H)=\{ \widetilde{X} \ |  \  X\in
\mbox{ind}(H-\mbox{mod}\vee H[1])\}$ [4].

Now we recall the representations of a species of a valued graph
from [8]. A valued graph $(\G,\b)$ is a finite set $\G$ (of
vertices) together with non-negative integers $d_{ij}$ for all
pair $i,j \in \G$ such that $d_{ii}=0$ and there exist positive
integers $\{\varepsilon _i\}_{i\in \G}$ satisfying
$$d_{ij}\varepsilon _j=d_{ji}\varepsilon _i,\  \ \mbox{for all }i,
\ j \in \G.$$
 A pair $\{i, j\}$ of vertices is called an edge of $(\G,\b)$ if
 $d_{ij}\not=0.$ An orientation $\Omega$ of a valued graph $(\G,\b)$
  is given by prescribing for each edge $\{i, j\}$ of $(\G,\b)$ an
  order (indicated by an arrow $i\rightarrow j$). Given an orientation
  $\Omega$
  and a vertex $k\in \G, $ we can define a new orientation
  $s_k\Omega$ of $(\G,\b)$ by reversing the direction of arrows
  along all edges containing $k$. A vertex $k\in \G$ is called
  a sink (or a source) with respect to $\Omega$ if there are no arrows
  starting (resp., ending) at vertex $k$.

  Let $K$ be a field and $(\G,\b, \Omega)$ a valued quiver.
  From now on, we shall always assume that
  $\G$ contains no cycles. Let
  $\m=(F_i, {}_iM_j)_{i,j\in \G}$ be a reduced $K-$species of
  type $(\G,\b, \Omega); $ that
  is, for all $i, j \in \G$, $_iM_j$ is an $F_i-F_j-$bimodule,
  where $F_i$ and $F_j$
   are division rings which are finite dimensional vector spaces over $K$
    and dim$(_{i}M_{j})_{F_j}=d_{ij}$
   and dim$_{K}F_i=
   \varepsilon_i$.
    A $K-$representation $V=(V_i,\varphi_{\alpha})$ of
   $(\m , \G, \Omega)$ consists of $F_i-$ vector space
   $V_{i}, i\in \G$, and of
   a $F_j-$linear map $_j\varphi_i: V_i\otimes {}_iM_j\rightarrow V_j$
   for each arrow $i\rightarrow j$. Such representation is called
   finite dimensional if $\sum _{i\in \Gamma}\mbox{dim}_{K}V_i<\infty.$ The
   category of finite-dimensional representations of $(\m , \G, \Omega)$ over
   $K$ is denoted by rep$(\m , \G, \Omega)$.

Now we fix a $K-$species $\m $ of type $(\G,\b, \Omega)$. In the
rest of the paper, we always speak of
   the valued quiver $(\m , \G, \Omega)$ instead of $(\G,\b, \Omega)$.
   Given a
sink, or a source $k$ of the valued quiver $(\m , \G, \Omega)$, we
are going to recall the Bernstein-Gelfand-Ponomarev reflection
functor (shortened as BGP reflection functor) $S^{\pm}_k$:

    $$S^+_k :\  \mbox{rep}(\m ,\G,  \Omega) \longrightarrow  \mbox{rep}(\m
    , \G, s_k\Omega), $$ respectively
$$S^-_k :\  \mbox{rep}(\m ,\G,  \Omega) \longrightarrow  \mbox{rep}(\m
,\G, s_k\Omega). $$

\medskip

 For any representation $V=(V_i, \phi _{\alpha})$ of
$(\m , \G, \Omega)$, the image of it under $S^+_k$ is by
definition, $S^+_kV=(W_i,\psi_{\alpha}),$ a representation of $(\m
, \G, s_k\Omega)$, where $W_i=V_i$ when $i\not=k;$ and $W_k$ is
the kernel in the diagram:
$$\begin{array}{lcccccccc}
(*)&&0&\longrightarrow& W_k&\s{({}_j\chi
_{k})_{j}}{\longrightarrow}& \oplus {}_{j\in \G} V_j\otimes
{}_jM_k &\s{({}_k\phi _{j})_j}{\longrightarrow}&V_k
\end{array}$$
$\psi _{\alpha }=\phi _{\alpha}$ when the ending vertex of
$\alpha$ is not $k$;  and when the ending vertex of $\alpha$ is
$k$, $\psi
 _{s_k\alpha }={}_j\bar{\chi}_{k}: W_k\otimes {}_kM_j\rightarrow X_j,$
 where ${}_j\bar{\chi}_{k}$ corresponds to ${}_j\chi _{k}$ under
 the isomorphism Hom$_{F_j}(W_{k}\otimes
 {}_kM_j,V_j)\approx\mbox{Hom}_{F_i}(W_k, V_j\otimes
 {}_jM_i ).$

If \textbf{f}$=(f _{i}): V\rightarrow V'$ is a morphism in rep$(\m
,\G,  \Omega)$, then $S^+_k($\textbf{f}$)=$\textbf{g}$=(g _i)$,
where $g _i =f _i $ for $i\not=k$ and $g _k: W_k\rightarrow W_k'$
as the restriction of $\oplus _{j\in \G}(f _j\otimes 1)$ given in
the following commutative diagram:

\[ \begin{CD}
0@>>>W_k @>(_j\chi_{k})_{j}>>\oplus_{j\in \G} V_j\otimes {}_jM_k @>(_k\phi_j)_j>>V_k\\
@VVV@VVg _k V @VV\oplus _j(f _j\otimes 1)V @VVf _k V \\
0@>>>W_k'@>(_j\chi_{k}')_{j}>>\oplus_{j\in \G} V_j'\otimes
{}_jM_k@>(_k\phi_j')_j>>V_k'
    \end{CD} \]

\medskip

If $k$ is a source, the definition of $S^-_kV$ is dual to that of
$S^+_kV$, we omit the details and refer to [8].

Let $k$ be a sink and $P_i$ the  indecomposable projective representation of $(\m,\G, \Omega)$
corresponding to vertex $i\in \G$. Let $T=\oplus _{i\in
\G-{k}}P_i\oplus \tau ^{-1}P_k$ and $H=\oplus _{i\in
\G}P_i$. Then $T$ is a tilting module in
rep$(\m,\G, \Omega)$ [1], $S^+_k=\mbox{Hom}(T,-)$. It induces an equivalence from add$T$ to add$H'$
where $H'$ is the tensor algebra of $(\m,\G, s_k\Omega)$, and
 induces a triangle equivalence $\mbox{Hom}_H(T,-): K^b(\mbox{addT})\rightarrow K^b(\mbox{addH}')$.  As in [12], the composition of functors
indicated as the following arrows:
$$\mathbf{K}^b(\mbox{add}T)\hookrightarrow
\mathbf{K}^b(H)\rightarrow D^b(H)$$ is a
triangle equivalence.  It is easy to see that $S^{+}_k$ and
$S^{-}_k$ commutes with the shift functor $[1]$. Since
$D^b(H)$ has Auslander-Reiten triangles and $S^+_k$ or
$S^{-}_k$ sends AR-triangles to AR-triangles (compare to Theorem
4.6 in Chapter I in [12]), $S^+_k$ and $S^{-}_k$ commute with $\tau
.$

\medskip

We summarize these facts in the following lemma.
\medskip

 \textbf{Lemma 2.1. }{\it Let $k$ be a sink (or a source) of a valued quiver
 $(\m , \G, \Omega).$
  Then $S^+_k$ (resp., $S^-_k$ )
   induces a triangle equivalence from
  $ D^{b}(H)$ to  $ D^{b}(H')$ which is denoted also by $S^+_k$ (resp. $S^-_k)$;
  and $S^{\pm}_k$ commutes with the shift functor
  $[1]$ and the AR-translation $\tau$.}
\medskip

In the following, we assume that the standard equivalence
$G: D^b(H)\rightarrow D^b(H)$ satisfies the
conditions $(g1)$ and $(g2)$. Then $G'=S^+_kGS^-_k$ is also a standard equivalence of
 $D^b(H')$ which satisfies $(g1)$ and $(g2)$.
 We define a functor $R(S^+_k)$ from
$D^{b}(H)/G$ to $D^{b}(H')/G'$ as follows: Let
$\widetilde{X}\in D^b(H)/G$ with $X\in
D^{b}(H)$. Let $X_{T}$ be one of the complexes in $C^b(\mbox{addT})$ which are quasi-isomorphic to $X$, where
 $C^b(\mbox{addT})$ denotes the category of complexes with finitely many non-zero components and all components belong to add$T$. We set
$R(S^+_k)(\widetilde{X})=\widetilde{S^+_k(X_{T})}.$ For morphism
$\tilde{f}:\ \widetilde{X}\rightarrow \widetilde{Y}$, we set
$R(S^+_k)(\tilde{f}):\widetilde{S^+_k(X_{T})}\rightarrow
\widetilde{S^+_k(Y_{T})}$ to be the map
$\widetilde{S^+_i(f_{T})},$ where $f_{T}$ is that one induced
from $f$ under the quasi-isomorphism from $X$ to $X_{T}$.
\medskip

 We prove that $R(S^+_k)$ is a triangle equivalence (compare Section 9.4 in [14]).
\medskip

 \textbf{Theorem 2.2. }{\it Let $k$ be a sink (or a source) of a
 valued quiver $(\m , \G, \Omega).$ Then $R(S^+_k)$ (resp., $R(S^-_k)$)
  is a triangle equivalence from $D^{b}(H)/G$ to
  $D^{b}(H')/G'.$}

\medskip

\textbf{Proof.} First of all, we verify the definition is
well-defined: For $\widetilde{X}= \widetilde{Y}\in
D^b(H)/G$ with $X, Y \in D^{b}(H),$ we have
that $Y=G^i(X)$ for some integer $i$. It follows that $Y_{T}=
G^i(X_{T})$ in $D^b(H).$ By applying $S^+_k$ to the
two complexes above, we have that $S^+_k(Y_{T})=S^+_kG^iS^-_k(S^+_k(X_{T}))=
G'^i(S^+_k(X_{T})).$ It follows that
$\widetilde{S^+_k(Y_{T})}=\widetilde{S^+_k(X_{T})},$ i.e.
$R(S^+_k)(\widetilde{X})= R(S^+_k)(\widetilde{Y})$. The action of
$R(S^+_k)$ on morphisms is induced by $S^+_k$ on morphisms in
$D^b(H)$ in the way indicated in the following
commutative diagram:

\[ \begin{CD}
 \oplus_{i\in \mathbf{Z}}\Hom _{D^b(H)}(G^i(X_{\mbox{T}}),
 Y_{\mbox{T}})
  @>S^+_k>>  \oplus_{i\in \mathbf{Z}}\Hom _{D^b(H')}
  (G'^i(S^+_k(X_{\mbox{T}})),
   S^+_k(Y_{\mbox{T}}))\\
@V\wr VV  @VV\wr V  \\
\Hom _{D^b(H)/G}(\widetilde{X},\widetilde{Y}) @>R(S^+_k)
>> \Hom _{D^b(H')/G'}(R(S^+_k)(\widetilde{X}),R(S^+_k)(\widetilde{Y}))
\end{CD} \]

 It is easy to verify that $R(S^+_k)$ and $R(S^-_k)$ satisfy:
 $R(S^+_k)\circ R(S^-_k)\approx \mbox{id}_{D^b(H')/G'}$ and
 $R(S^-_k)\circ R(S^+_k)\approx \mbox{id}_{D^b(H)/G}$.
 These show that $R(S^+_k)$ and $R(S^-_k)$ are equivalences.
 Now by using the result in section 9.4 of [14], we have that $R(S^+_k)$ sends triangles in $D^b(H)/G$ to
 triangles in
 $D^b(H')/G' $.  Therefore
  $R(S^+_k)$ is a triangle equivalence. The proof is finished.

 \medskip

When $G=\tau ^{-1}[1]$, we have the triangle equivalence
$R(S^+_k)$
 from the cluster category $\mathcal{C}(\Omega) =D^{b}(H)/G$ to
$\mathcal{C}(s_k\Omega) =D^{b}(H')/G.$ And when $G=[2],$ we have the triangle equivalence from
the root category $D^{b}(H)/[2]$ to the root category $D^{b}(H')/[2]$ (compare [21]). \medskip

 Let
$P_i$ (or $P_i'$) be the indecomposable projective representations
in $H-$mod (resp. $H'-$mod) corresponding to the
vertex $i\in \G_0$, $E_j$ (or $E_j'$) the simple
 $H-$module (resp. simple $H'-$module)
 corresponding to the vertex $j$.

\medskip

 \textbf{Corollary 2.3. }{\it Let $k$
be a sink of a
 valued quiver $(\m , \G, \Omega).$
 Then $R(S^+_k)$ is a triangle equivalence from
  $\mathcal{C}(\Omega)$ to
  $\mathcal{C}(s_k\Omega).$ Moreover for $X\in \mbox{ind}H,$
$ R(S^+_k)(\widetilde{X})= \left\{ \begin{array}{ll}
\widetilde{ P_k'}[1] &\mbox{ if  }X\cong E_k  \\
\widetilde{ S^+_k(X)} & \mbox{otherwise; }
\end{array}\right.$
  and for $j\neq k$, $R(S^+_k)(\widetilde{P_j}[1])=\widetilde{P_j'}[1],$
    and $R(S^+_k)(\widetilde{P_k}[1])=\widetilde{E_k'}.$
     }
\medskip

\textbf{Proof.} From Theorem 2.2., $R(S^+_k)$ is a triangle equivalence from the cluster category $\mathcal{C}(\Omega)$ to
$\mathcal{C}'(s_k\Omega)$.  Now we
prove that $R(S^+_k)(\widetilde{E_k})=
  \widetilde{P_k'}[1].$ Since $k$ is sink, we have AR-sequence
  $ (*):\ 0\rightarrow
E_k\rightarrow
  X\rightarrow\tau^{-1}E_k\rightarrow 0$ in $H-$mod with
   $X$ and $\tau^{-1}E_k$ being in add$T$ [1]. Since $S^+_k$ is
 a left exact functor, we have the exact sequence
 $0\rightarrow
  S^+_k(X)\rightarrow S^+_k(\tau^{-1}E_k)$ in $H'-$mod, in which
  the cokernel of the injective map is $E_k'.$  As the stalk complex
   of degree $0$, $E_k^{\bullet}$ is isomorphic to the complex:
  $\cdots\rightarrow 0\rightarrow X\rightarrow \tau^{-1}E_k\rightarrow 0
  \rightarrow\cdots$
  in $D^b(H).$ By applying $S^+_k$ to the complex above,
   we have
  that $S^+_k(E_k^{\bullet})=\cdots\rightarrow 0\rightarrow S^+_k(X)
  \rightarrow S^+_k( \tau^{-1}E_k)\rightarrow0\rightarrow\cdots.$
   It follows that the complex $\cdots\rightarrow 0\rightarrow
S^+_k(X)\rightarrow S^+_k(\tau ^{-1}E_k)\rightarrow 0
\rightarrow\cdots$ is quasi-isomorphic to the stalk complex
$E_k'^{\bullet}[-1]$  of degree $-1$. It follows
 $R(S^+_k)(\widetilde{E_k})=\widetilde{E_k'}[-1].$ Since $\tau
\widetilde{P_k'} =\widetilde{E_k'}[-1], $
$R(S^+_k)(\widetilde{E_k})=\tau \widetilde{P_k'} =
\widetilde{(\tau [-1])(P_k')}[1]=\widetilde{P_k'}[1].$
  In the derived category $D^b(H)$,
  we have that $S^+_k(P_i)=P_i'$ for any $i\neq k$,  $S^+_k(E_k[1])=E_k'.$
      It
follows that $R(S^+_k)(\widetilde{P_i})=\widetilde{P_i'}$ for any
$i\neq k$ and $R(S^+_k)(\widetilde{P_k}[1])=\widetilde{E_k'}$. The
proof is finished.

\medskip

\textbf{Remark 2.4.} We leave the dual statement for a source $k$
 to the reader. \medskip

 \textbf{Definition 2.5. } When $k$ is a sink (or a source) of
 $(\m, \G, \Omega)$, the
functor $R(S^+_k)$ (resp. $R(S^-_k)$) in Corollary 2.3. is called a BGP-reflection functor in the cluster category
 $\mathcal{C}(\Omega)$.
\medskip

Let $\G$ be a classical Dynkin quiver, i.e. one of the types
$ADE$. Then the automorphisms of $D^b(H)$ are of the form:
$[n],\ \tau ^n, $ or of the
  form $[n]\tau ^{m}$ for any $m, n \in \mathbf{Z}$ (compare
  [20])

\medskip

 \textbf{Corollary 2.6. }{\it Let $\G$ be a Dynkin quiver and $k$
 a sink (or a source) of it. Then for any automorphism $G$  of
  $D^b(H)$ which is not of the forms
 $([1]\tau ^{m})^t$, where $t$ is an integer and $m= (n+1)/2$ if the underlying diagram of
 $\G $ is of type
  $A_n$;  or $m= 6$ if the underlying diagram of
 $\G $ is of type $E_6$, $R(S^+_k)$ (resp., $R(S^-_k)$)
  can be defined and it is a triangle equivalence from
  $D^{b}(H)/G$ to $D^{b}(H')/G'.$}
  \medskip

  \textbf{Proof.} It follows from Proposition 3.3.2 in [20] that
  any automorphism $G$
   of $D^b(H)$ is generated by $\tau$ and $[1]$.
  Therefore $G$ commutes with $S^+_k.$  For automorphism $G$ indicated
  in the corollary, $G$ satisfies the conditions
  $(g1)$ and $(g2)$, hence the orbit category
  $D^b(H)/G$ exists [14]. Then by Theorem 2.2., $R(S^+_k)$
  exists and is a triangle equivalence. The proof is finished.
  \medskip


 \begin{center}

\textbf{3. Applications to cluster combinatorics.}
\end{center}

\medskip

In this section, we always assume that $H$ is the tensor algebra of a
valued quiver $(\m,\G, \Omega)$ over a field $K,$  with underlying
graph $\G,$ where $\G$ is not necessarily connected. We denote by
$\mathcal{A}=\mathcal{A}(\G)$ the corresponding cluster algebra
when $\G$ is of Dynkin type (simply-laced or non-simply-laced),
by $\Phi$ the set of roots of the corresponding Lie algebra, and
 $\Phi_{\geq -1}$ the set of almost positive roots, i.e. the
positive roots together with the negatives of the simple roots.
The elements of $\Phi_{\geq -1}$ are in 1--1 correspondence with
cluster variables of $\mathcal{A}$ (Theorem 1.9. [10]), such 1--1
correspondence is denoted by $\mathcal{P}$. Fomin and Zelevinsky
[11] associate a nonnegative integer $(\alpha||\beta)$, known as
the { compatibility degree}, to each pair $\alpha,\beta$ of almost
positive roots. This is defined in the following way. Let $s_i$ be
the Coxeter generator of the Weyl group of $\Phi$ corresponding to
$i$, and let $\sigma_i$ be the permutation of $\Phi_{\geq -1}$
defined as follows:
$$(3.1)  \ \ \ \ \sigma_i(\alpha)=
\left\{ \begin{array}{ll} \alpha & \alpha=-\alpha_j,
\ j\not=i \\ s_i(\alpha) & \mbox{otherwise.} \end{array}\right.$$

The $\sigma_i$'s are called "truncated simple reflections" of
$\Phi_{\geq -1}$. They are one of the main ingredients of
constructions in [11] (see also [15]).  Let $\G =\G ^+\sqcup \G
^-$ be a partition of the set of vertices of $\G$ into completely
disconnected subsets and define:
$$(3.2)\ \ \ \ \tau_{\pm}=\prod_{i\in \G ^{\pm}}\sigma_i.$$
 Denote by $[\beta:\alpha _i]$ the coefficient of $\alpha_i$ in the expression of $\beta$ in simple roots $\alpha _1,
\cdots , \alpha _n$.
Then $(\ ||\ )$ is uniquely defined by the following two
properties:

 $$\begin{array}{lllcl}(3.3)& &(-\alpha_i||\beta)&=&\mbox{max}([\beta:\alpha _i],
 0),\\
  (3.4)&&(\tau_{\pm}\alpha||\tau_{\pm}\beta )&=&(\alpha ||\beta),\end{array}$$
  for any $\alpha , \beta \in \Phi_{\geq -1},$ any $i\in \G $.

 A pair $\alpha \ \beta$ in $\Phi_{\geq -1}$ are called compatible if
$(\alpha ||\beta)=0$. Associated to the finite root system $\Phi$,
the simplicial complex $\Delta (\Phi)$ is defined in [11].
$\Delta (\Phi)$ has $\Phi_{\geq -1}$ as the set of vertices, its
simplices are mutually compatible subsets of $\Phi_{\geq -1}.$ The
maximal simplices of $\Delta (\Phi)$ are called the clusters
associated to $\Phi.$ This simplicial complex $\Delta (\Phi)$ is
called generalized associahedron (compare [5, 6, 10, 11]).

 In this
section, we will first show that the truncated simple reflections
$\sigma _i$ on $\Phi_{\geq -1}$ can be realized by the BGP-reflection
functors $R(S^+_i)$ in the corresponding cluster category. Then, by using these BGP-reflection functors
and combining tilting theory in cluster categories developed in
[4], we give a unified quiver-interpretation of certain
combinatorics about clusters associated to arbitrary root systems
of (simply-laced or non-simply-laced) semisimple Lie algebras in
[11]. This extends, in a different way, the quiver-theoretic
interpretation of certain combinatorics about clusters in the
simply-laced case given by Marsh, Reineke
 and Zelevinsky in [15]. They use decorated representations.

Let $\{ e_i\ | \ i\in \G \}$ be a complete set of primitive
idempotents of a hereditary algebra $H$. For any subgraph $J$ of
$\G$, we set $I=HeH$ the hereditary ideal of H, where $e=\sum _{i
\in \G -J}e_i$. Then quotient algebra $A=H/I$ has a complete set
of primitive idempotents $\bar{e_i}$ $i\in J$. $A-$mod is a full
subcategory of $H-$mod consisting of $H-$modules annihilated
 by $I$ or in other words, consisting of
  $H-$modules whose composition factors are
 $E_i$ with $i\in J$. It follows from [8, 16] that
 Ext$^i_A(X,Y)=\mbox{Ext}^i_{H}(X,Y)$ for any $X, Y \in
 A-\mbox{mod}$ and any $i$. It follows that $A$ is also a hereditary
 algebra which is Morita equivalent to the tensor algebra
 of $(\m |_J ,J, \Omega |_J).$ These facts are summarized in the following
 proposition.
 \medskip

\textbf{Proposition 3.1.} {\it $D^b(A)$ is a triangulated
subcategory of $D^b(H)$ and $\mathcal{C}(A)$ is a triangulated
subcategory of $\mathcal{C}(H)$.}
\medskip

\textbf{Proof.} $A$ is hereditary and
Ext$^i_A(X,Y)=\mbox{Ext}^i_{H}(X,Y)$ for any $X, Y \in
 \mbox{mod}A$ and any $i$. This gives us that $D^b(A)\subseteq
 D^b(H)$ is a full triangulated subcategory of $D^b(H)$. It follows that
 the cluster category $\mathcal{C}(A)$ is a full triangulated
 subcategory of $\mathcal{C}(H)$.
  The proof is finished.
\medskip

 We recall the notation of exceptional sets and of tilting sets in
$\mathcal{C}(\Omega)$ in [4].
 A subset $B$ of ind$\mathcal{C}(\Omega)$ is called { exceptional } if Ext$^1_{\mathcal{C}(\Omega)}(X,Y)=0$ for any $X,\ Y\in
  B.$ A maximal exceptional set is called a tilting set. A subset of
  $\mathcal{C}(\Omega)$ is a tilting set if and only the direct
  sum of all objects in $B$ is a basic tilting object [4]. Then any tilting set
  contains exactly $|\G|$ many objects. One can associate to $\mathcal{C}(\Omega)$ a simplicial complex
 $\Delta(\Omega)$ as follows: $\Delta(\Omega)$ has ind$\mathcal{C}(\Omega)$ as the set of vertices, its simplices are the
exceptional sets in ind$\mathcal{C}(\Omega)$. It is easy to see that its maximal simplices are exactly tilting sets [4].  One can
 also associate to $\mathcal{C}(\Omega)$ a tilting graph $\Delta _{\Omega}$
   whose
  vertices are the basic tilting objects, and where there is an
  edge between two vertices if the corresponding tilting objects
  have all but one indecomposable summands in common. Tilting
  graphs associated to a hereditary algebra were studied by
  C.Riedtmann and A.Schofield [17] and L.Unger [19],
  also Happel, Unger [13].

In general, BGP-reflection functors preserve exceptional sets and tilting sets.
 \medskip

 \textbf{Proposition 3.2.}{ \it Let $k$
be a sink (or a source) of a valued quiver $(\m , \G, \Omega)$ of
any type. Then the BGP- reflection functor $R(S^+_k)$ ($R(S^{-}_k)$,
resp.) gives a 1-1 correspondence from the set of exceptional
sets in
 ind$\mathcal{C}(\Omega)$ to that in ind$\mathcal{C}(s_k\Omega)$, under
 this correspondence, tilting sets go to tilting sets. In particular if
$(\m , \G, \Omega)$ and $(\m , \G, \Omega ')$ are two valued
quivers of the same type $\G$, then the simplicial complexes $\Delta(\Omega)$ and $\Delta(\Omega')$ are isomorphic and the tilting graphs $\Delta
_{\Omega}$ and $\Delta _{\Omega'}$ are isomorphic.}
\medskip

\textbf{Proof.} Suppose $k$ is a sink. Since $R(S^+_k)$ and
$R(S^-_k)$ are inverse equivalences between $\mathcal{C}(\Omega)$
and $\mathcal{C}(s_k\Omega)$,  $$ \Ext
^1_{\mathcal{C}(\Omega)}(X,Y)=\Ext
^1_{\mathcal{C}(s_k\Omega)}(R(S^+_k)(X),R(S^+_k)(Y)), \mbox{ for
any } X,\  Y \in \mathcal{C}(\Omega).$$ It follows that
  $R(S^+_k)$ and $R(S^-_k)$ give inverse maps between the sets of
  exceptional sets in
 ind$\mathcal{C}(\Omega)$ and in ind$\mathcal{C}(s_k\Omega).$
  A exceptional set is a tilting set if and only so is
   its image under $R(S^+_k)$.  For any two valued quivers with the same
graph, one can get an admissible sequence $i_1, \cdots , i_n$ such
that $\Omega' =s_{i_n}\cdots s_{i_1}\Omega$ with $i_k$ is the sink
of $s_{i_{k-1}}\cdots s_{i_1}\Omega .$ For each {k}, we have that
the fact of equivalence of $R(S^+_{i_k})$ implies  $\Delta
_{s_{i_{k-1}}\cdots s_{i_1}\Omega}\simeq \Delta
_{s_{i_k}s_{i_{k-1}}\cdots s_{i_1}\Omega}$ and $\Delta
(s_{i_{k-1}}\cdots s_{i_1}\Omega)\simeq \Delta
(s_{i_k}s_{i_{k-1}}\cdots s_{i_1}\Omega)$. Therefore $\Delta
_{\Omega}\simeq \Delta _{\Omega'}$ and  $\Delta
(\Omega)\simeq \Delta (\Omega')$. The proof is finished.

\medskip

  Now we recall the decorated quiver representations from [15].
   Let $Q$ be a Dynkin quiver with vertices $Q_0$ and arrows $Q_1$.
   The ``decorated'' quiver $\widetilde{Q}$ is the quiver $Q$
   with an extra copy
$Q_0^-=\{i_-\,:\,i\in Q_0\}$ of the vertices of $Q$ (with no
arrows incident with the new copy). A module $M$ over
$k\widetilde{Q}$ can be written in the form $M^+\oplus V$, where
$M^+=\oplus_{i\in Q_0}M^+_i$ is a $KQ$-module, and $V=\oplus_{i\in
Q_0}V_i$ is a $Q_0$-graded vector space over $K$. Its signed
dimension vector, $\mathbf{sdim}(M)$ is the element of the root
lattice of the Lie algebra of type $Q$ given by
$$\mathbf{sdim}(M)=\sum_{i\in Q_0}\dim(M^+_i)\alpha_i-
\sum_{i\in Q_0}\dim(V_i)\alpha_i,$$ where
$\alpha_1,\alpha_2,\ldots ,\alpha_n$ are the simple roots. By
Gabriel's Theorem, the indecomposable objects of
$K\widetilde{Q}$-mod are parameterized, via $\mathbf{sdim}$, by
the almost positive roots, $\Phi_{\geq -1}$, of the corresponding
Lie algebra. The positive roots correspond to the indecomposable
$KQ$-modules, and the negative simple roots correspond to the
simple modules associated with the new vertices. We denote the
simple module corresponding to the vertex $i_-$ by $E_i^-$. Let
$M=M^+\oplus V$ and $N=N^+\oplus W$ be two
$K\widetilde{Q}$-modules. The symmetrized $\Ext^1$-group for this
pair of modules is defined to be:
\begin{eqnarray*}
E_{KQ}(M,N) & := & \Ext^1_{KQ}(M^+,N^+)\oplus \Ext^1_{KQ}(N^+,M^+)
\oplus \\
& & \Hom^{Q_0}(M^+,W) \oplus \Hom^{Q_0}(V,N^+),
\end{eqnarray*}
where $\Hom^{Q_0}$ denotes homomorphisms of $Q_0$-graded vector
spaces.

The map $\psi_Q$ from $\ind\mathcal{C}(KQ)$ to the set of
isomorphism classes of indecomposable $K\widetilde{Q}$-modules  is
defined in [4] as follows. Let
$\widetilde{X}\in\ind\mathcal{C}(KQ)$. It can be assumed that one
of the following cases holds:
\begin{enumerate}
\item $X$ is an indecomposable $KQ$-module $M^+$.
\item $X=P_i[1]$ where $P_i$ is the indecomposable projective $KQ$-module
corresponding to vertex $i\in Q_0$.
\end{enumerate}
We define $\psi_Q(\widetilde{X})$ to be $M^+$ in Case (1), and to
be $E_i^-$ in Case (2).

Then the map $\psi_Q$ is a bijection between $\ind\mathcal{C}(KQ)$
and the set of isomorphism classes of indecomposable
$K\widetilde{Q}$-modules (i.e.\ indecomposable decorated
representations). If we denote by
 $\gamma_Q:=\mathbf{sdim}\circ \psi_Q$, then it is a bijection between
$\ind\mathcal{C}(KQ)$ and $\Phi_{\geq -1}$ (and thus induces a
bijection between $\ind\mathcal{C}(KQ)$ and the set of cluster
variables). For $\alpha\in \Phi_{\geq -1}$ we denote by
$M_Q(\alpha)$ the element of $\ind\mathcal{C}(KQ)$ such that
$\gamma_Q(M_Q(\alpha))=\alpha$. It was proved in [4] that
$$E_{KQ}(\psi_Q(\widetilde{X}),\psi_Q(\widetilde{Y}))\simeq
\Ext^1_{\mathcal{C}(KQ)}(\widetilde{X},\widetilde{Y}), \ \
\mbox{for } X,Y \in \mathcal{D}.$$
\medskip

Now we return to the general case. Let $(\m , \G, \Omega)$ be a
Dynkin valued quiver. We extend first the bijection $\gamma _{Q}$
to the general case $\gamma_{(\m, \G,\Omega)}$(which is denoted
for simplicity by $\gamma_{\Omega}$) from $\ind\mathcal{C}$ to
$\Phi_{\geq -1}$ by defining: Let $X\in \mbox{ind}(H-\mbox{mod}\vee
H[1]).$
$$\gamma_{\Omega}(\widetilde{X})=\{
\begin{array}{lrl}\mathbf{dim}X & \mbox{ if } & X\in \mbox{ind}H;\\
&&\\
-\mathbf{dim}E_i& \mbox{ if } &X=P_i[1],\end{array}$$ where
$\mathbf{dim}X$ denotes the dimension vector of $H-$module $X$. It
is easy to see the map $\gamma_{\Omega}$ is a bijection and it is
dependent on the orientation $\Omega$ of $\G$.
\medskip

Let $\alpha ,\beta \in\Phi_{\geq -1}$ and  $M(\alpha)$, $M(\beta)$
the indecomposable objects in $\mathcal{C}(\Omega)$ corresponding
to $\alpha ,\beta$ under the bijection $\gamma_{\Omega}$. For any
pair of objects $M, \ N $ in $\mathcal{C}(\Omega)$, Hom$_{\mathcal{C}(\Omega)}(M,N)$
is a left End$_{ \mathcal{C}(\Omega)}M-$module
(here the composition of maps $f:X\rightarrow Y$ and $g:Y\rightarrow Z$ is $f\circ g: X\rightarrow Z$). Therefore
Ext$^1_{\mathcal{C}(\Omega)}(M,N)$ is a left End$_{
\mathcal{C}(\Omega)}M-$module. For an algebra $A$ and an $A-$module $X$, we denote by $l(_AX)$ the length of $A-$module $X$.
 It is easy to see that $R(S^+_k)$ induces an isomorphism from
  End$_{ \mathcal{C}(\Omega)}M$ to End$_{ \mathcal{C}(s_k\Omega)}(R(S^+_k)M)$.
  Under this isomorphism, $R(S^+_k)$ induces an End$_{ \mathcal{C}(\Omega)}M$
  ($\cong \mbox{ End}_{ \mathcal{C}(s_k\Omega)}(R(S^+_k)M))$
-module
  isomorphism
 between Ext$^1_{\mathcal{C}(\Omega)}(M,N)$ and
  Ext$^1_{\mathcal{C}(s_k\Omega)}(R(S^+_k)(M),R(S^+_k)(N)),$ for any
  sink $k$. Similar isomorphisms
 hold if $k$ is a source.
\medskip

 \textbf{Definition 3.3.} For any two almost positive roots $\alpha ,
 \beta \in \Phi _{\geq -1}$, we define the $\Omega-$ compatibility
 degree $(\alpha||\beta)_{\Omega}$ of $\alpha , \beta$ by
  $$(\alpha||\beta)_{\Omega}=l(_{\mbox{End}M(\alpha)}
  \mbox{Ext}^1_{\mathcal{C}
  (\Omega)}(M(\alpha), M(\beta))).$$
\medskip

Note that if $(\G, \Omega)$ is a simply-laced Dynkin quiver, then
the $\Omega-$ compatibility
 degree $(\alpha||\beta)_{\Omega}$ of $\alpha , \beta$ equals
$ \dim_K  \mbox{Ext}^1_{\mathcal{C}
  (\Omega)}(M(\alpha), M(\beta)).$
  \medskip

 We now prove the first main result of the paper.
\medskip

 \textbf{Theorem 3.4. }{\it Let $(\m , \G, \Omega)$ be a valued
 Dynkin quiver and
 $k$ a sink (or a source). Then we have the commutative diagram:
\[ \begin{CD}
\Phi_{\geq -1} @>\sigma _k>> \Phi_{\geq -1}\\
@V\gamma ^{-1}_{\Omega} VV  @VV\gamma ^{-1} _{s_k\Omega}V  \\
\mbox{ind}\mathcal{C}(\Omega) @>R(S^+_k)
>(R(S^-_k),resp.)> \mbox{ind}\mathcal{C}(s_k\Omega)
\end{CD} \]

Moreover
$(\alpha||\beta)_{\Omega}=(\sigma_k(\alpha)||\sigma_k(\beta))_{s_k\Omega}.$}

\medskip

\textbf{Proof.}
 Let $\alpha \in
\Phi_{\geq -1}$ be a positive root. Then  $\sigma _k(\alpha)=-\alpha _k$ when $\alpha =\alpha _k$, and $\sigma _k(\alpha)=s_k(\alpha)$ when $\alpha$
 is a positive root other than $\alpha _k$. It follows that $\gamma ^{-1} _{s_k\Omega}\sigma
_k(\alpha)$ is $\widetilde{P_k'}[1]$ or $\widetilde{S^+_k(X)},$
respectively, where  $X$ is the unique indecomposable
representation with $\mathbf{dim}X=\alpha$ which does exist by
Gabriel's theorem [7]. On the other side, $R(S^+_k)\gamma ^{-1}
_{\Omega}(\alpha)$ equals to $R(S^+_k)(\widetilde{E_k})$ or
$R(S^+_k)(\widetilde{X})$ according to $\alpha $ is simple root
$\alpha _k$ or not. Then it follows from Corollary 2.3 that
$\gamma ^{-1} _{s_k\Omega}\sigma _k(\alpha)=R(S^+_k)\gamma
^{-1}_{\Omega}(\alpha).$ We now prove the equality above for
$\alpha $ a negative root. Let $\alpha=-\alpha_i$ for $i\in \G .$
 Then we have that
$$\gamma ^{-1} _{s_k\Omega}\sigma
_k(-\alpha_i)=\{ \begin{array}{crc}\widetilde{E_k'}& \mbox{ if }&
i=k\\
\widetilde{P_i'[1]}& \mbox{  if }& i\not= k.\end{array}$$ Again
from Corollary 2.4, we have that $R(S^+_k)\gamma ^{-1}
_{\Omega}(-\alpha
_k)=R(S^+_k)(\widetilde{P_k[1]})=\widetilde{E_k'[1]}$ and for
$i\not= k$, $R(S^+_k)\gamma ^{-1} _{\Omega}(-\alpha
_i)=R(S^+_k)(\widetilde{P_i[1]})=\widetilde{P_i'[1]}.$
 This finishes the proof of the commutativity of the diagram.
   By definition,  $(\sigma_k(\alpha)||\sigma_k(\beta))_{s_k\Omega}
   \newline
= l(_{\mbox{End}M(\sigma_k(\alpha))}
  \mbox{Ext}^1_{\mathcal{C}
  (s_k\Omega)}(M(\sigma_k(\alpha)),
 M(\sigma_k(\beta)))).$ On the other hand, it follows from the
  commutative diagram which is proved above,
 that $$(\sigma_k(\alpha)||\sigma_k(\beta))_{s_k\Omega} =
l(_{\mbox{End}R(S^+_k)(M(\alpha))}
  \mbox{Ext}^1_{\mathcal{C}
  (s_k\Omega)}(R(S^+_k)(M(\alpha)),R(S^+_k)(M(\beta)))).$$
 The right hand of the equality equals
  $l(_{\mbox{End}M(\alpha)}
  \mbox{Ext}^1_{\mathcal{C}
  (\Omega)}(M(\alpha), M(\beta)))$ since $R(S^+_k)$
 is a triangle equivalence. Therefore
 $(\alpha||\beta)_{\Omega} =(\sigma_k(\alpha)||\sigma_k(\beta))_{s_k\Omega}.$
 The proof is finished.
\medskip

\textbf{Remark 3.5.} If the valued quiver $(\m , \G, \Omega)$ is
simply-laced, then we have the following commutative diagram:

\[ \begin{CD}
\mbox{ind}\mathcal{C}(\G) @>R(S^+_k)
>(R(S^-_k),resp.)> \mbox{ind}\mathcal{C}(s_k\G)\\
@V\Psi_{\G} VV  @VV\Psi _{s_k\G}V\\
\mbox{indrep}\widetilde{\G }@>\Sigma ^+ _k>(\Sigma ^- _k,resp.)>
 \mbox{indrep}\widetilde{s_{k}\G} \\
@V\mathbf{sdim} VV  @VV\mathbf{sdim}V\\
\Phi_{\geq -1} @>\sigma _k>> \Phi_{\geq -1},
\end{CD} \]

 and the $\Omega$-compatibility degree of $\alpha$ and $\beta$
defined above is the same as defined in [15] (compare [4]).
This implies Theorem 4.7. there. We remark that the functors
$\Sigma ^+_k$ and $\Sigma ^-_k$
 defined in [15] are not equivalences.

 The next result shows that
 the $\Omega-$compatibility degree function on $\Phi
_{\ge -1}$ is independent of the orientation $\Omega$ of $\G$. It
is the same as that defined in [11] on $\Phi _{\ge -1}$. This
gives a unified form of compatibility degree in the language of
quiver representations.

\medskip

\textbf{Theorem 3.6.}{ \it  Let $(\m , \G, \Omega _0)$ be an alternative valued quiver of Dynkin type. The $\Omega _0-$compatibility degree function on $\Phi
_{\ge -1}$ is the same as the compatibility degree function given by
Fomin-Zelevinsky in [10, 11].}
\medskip

\textbf{Proof.} We have to verify the $\Omega _0-$compatible degree
function satisfies conditions $(3.3), (3.4)$. For any two
orientations on the same graph $\G$, one can get an admissible
sequence $i_1, \cdots , i_n$ such that $\Omega' =s_{i_n}\cdots
s_{i_1}\Omega,$ where $i_k$ is the sink of $s_{i_{k-1}}\cdots
s_{i_1}\Omega .$ Then, from Theorem 3.4., we have
$(\alpha||\beta)_{\Omega}=(\sigma_{i_n}\cdots \sigma_{i_1}(
\alpha)||\sigma_{i_n}\cdots \sigma_{i_1}(\beta))_{\Omega'}$. It
follows that $(\tau_{\varepsilon}\alpha||\tau_{\varepsilon}\beta
)_{\tau_{\varepsilon}(\Omega)}=(\alpha ||\beta)_{\Omega},$
  for any $\alpha , \beta \in \Phi_{\geq -1},$ any
  $\varepsilon\in \{-1, 1\}$. This proves that $(3.4)$ holds.
   Let $\beta \in \Phi_{\geq -1}$. Then
  $(-\alpha _i||\beta)_{\Omega _0}=l(_{\mbox{End}\widetilde{P_i[1]}}
  \mbox{Ext}^1_{\mathcal{C}
  (\Omega _0)}(\widetilde{P_i[1]},\widetilde{ M(\beta)}))=l(_{\mbox{End}\widetilde{P_i}}
  \mbox{Hom}_{\mathcal{C}
  (\Omega _0)}(\widetilde{P_i}, \widetilde{M(\beta)})).$ It equals $ l(_{\mbox{End}_HP_i}\mbox{Hom}_{H}(P_i, M(\beta)))$ (this follows from
  Proposition 1.7. in [4]) and then it equals $[\beta: \alpha_i]$
  if $\beta $ is a positive root, or $0$ otherwise. This proves that $(3.3)$ holds. The proof is
  finished.
\medskip

This theorem extends Proposition 4.2 in [4] since in the
simply-laced case, the $\G-$compatible degree defined in [15] is
also the same as the compatible degree [11]. Since Theorem
3.6, we denote $(\alpha||\beta)_{\Omega _0}$ just by
$(\alpha||\beta)$.
\medskip

  \textbf{Definition 3.7.} A subset $C$ of $\Phi _{\geq -1}$ is called
  compatible if $(\alpha||\beta)=0$ for all
 $\alpha , \ \beta \in C.$ The subset $C$ is called a cluster if
 it is a maximal compatible.
\medskip

\textbf{Definition 3.8.} The negative support $S(C)$ of a subset
$C$ of $\Phi _{\geq -1}$ is defined by  $S(C)=\{ i\in \G: -\alpha
_i \in C \}.$ The subset $C$ is called positive if $C\subset\Phi
_{> 0},$ i.e. $S(C)=\emptyset$ [11, 15].

\medskip

Combining Theorem 3.6. with Proposition 1.7 in [4], we reprove
 Propositions 3.3. 3.5. 3.6. in [11] in the language of
 quiver representations.
\medskip

 \textbf{Proposition 3.9.} {\it Let $\G$ be any Dynkin diagram, $k$
 a vertex of $\G$ and $\alpha, \ \beta$ almost positive roots.
 Then

 (1) $(\alpha||\beta)=(\beta||\alpha)$ if $\G$ is  a simply-laced
 Dynkin quiver;

(2) $\sigma _k$ sends a compatible subset to a compatible subset.
In particular, it sends clusters to clusters;

(3) If $\alpha $ and $\beta $ belong to $\Phi(J)_{\ge -1}$ for
some proper subset $J\subset \G,$  then their compatibility degree
with respect to the root subsystem $\Phi(J)$ is equal to
$(\alpha||\beta)$

(4) If $\G_1, \cdots, \G_r \subset \G$ are the connected
components of the coxeter graph, then the compatible subsets
(resp., clusters) for $\Phi(\G)_{\ge -1}$ are the disjoint unions
$A_1\coprod \cdots \coprod A_r$, where each $A_k$ is a compatible
subset (resp., clusters) for $\Phi(\G_k)_{\ge -1}.$

(5) For every subset $J\subset \G,$ the correspondence $C\
\mapsto C-\{-\alpha _i: i\in J\}$ is a bijection between the set
of all compatible subsets (resp., clusters) for $\Phi(\G)_{\ge
-1}$ with negative support $J$ and the set of all positive
compatible subsets (resp., clusters) for $\Phi(\G-J)_{\ge -1}$.}
\medskip

\textbf{Proof.} Let $(\m , \G, \Omega)$ be a valued quiver with
the underlying diagram $\G$ such that $k$ is a sink. Statement (1) follows from Ext$^1_{\mathcal{C}(\Omega)}(M(\alpha),
M(\beta))= \mbox{Ext}^1 _{\mathcal{C}(\Omega)}(M(\beta),M(\alpha))$
[4]. For statement (2), we note that
$(\sigma_k(\alpha)||\sigma_k(\beta))=(\sigma_k(\alpha)||\sigma_k(\beta))
_{s_k\Omega}=0$ if and only
 if $(\alpha||\beta )
_{\Omega}=(\alpha||\beta)=0.$  Statement (3) follows from
Proposition 3.1.  Statement (4) follows from (3),
Proposition 3.1. and the obvious fact
ind$\mathcal{C}(\G)=\mbox{ind}\mathcal{C}(\G _1)\coprod \cdots
 \coprod \mbox{ind}\mathcal{C}(\G _r).$  For the proof of (5), we assume
 $C=\{\tilde{X} : X\in \mbox{ind}H\}
 \coprod \{\tilde{P_{i}}[1] : i\in J \} $ is a compatible subset (resp., cluster) for $\Phi(\G)_{\ge
-1}$ with negative support $J$. Then
Ext$^1_{\mathcal{C}(H)}(\tilde{P_{i}}[1],\tilde{X})=0$. It follows
that $$\begin{array}{ll}\mbox{Hom}_H(P_i,
X)=\mbox{Hom}_{\mathcal{C}(H)}(\tilde{P_{i}},\tilde{X}) = &\\
 \mbox{Hom}_{\mathcal{C}(H)}(\tilde{P_{i}}[1],\tilde{X}[1]) =
 \mbox{Ext}_{\mathcal{C}(H)}^1(\tilde{P_{i}}[1],\tilde{X})& \\ =0.\end{array}$$
 Then $X\in \mbox{ind} A$ where $A=H/I$ is the quotient algebra of $H$
 whose modules are exactly the $H-$modules without composition factors $E_i$, with
$i\in J$ (compare Proposition 3.1.). Then $C-\{\tilde{P_{i}}[1] : i\in J \}
 = \{\tilde{X} : X\in
 \mbox{ind}H\}$ is a compatible subset (resp., cluster) of
  $\Phi(\G-J)_{\ge -1}$. Conversely,
  given a compatible subset (resp., cluster) $C_1=\{\tilde{X} : X\in
 \mbox{ind}A \}$ of
  $\Phi(\G-J)_{\ge -1}$, $C_1\coprod \{\tilde{P_{i}}[1] : i\in J \}$
  is a compatible subset (resp., cluster) for $\Phi(\G)_{\ge
-1}$ with negative support $J$. The proof is finished.

\medskip

As a consequence of Theorems 3.4., 3.6, we have the second main
result of the paper which is a generalization of Theorem 4.5.
in [4] and confirm positively the Conjecture 9.1. in all
Dynkin types.
\medskip

\textbf{Theorem 3.10. }{\it Let $(\G,\b, \Omega)$ be any Dynkin
valued quiver, $\Phi _{\geq -1}$ the  set of almost positive roots
of the corresponding Lie algebra. Then the bijection $\gamma
_{\Omega} : \ind\mathcal{C}(\Omega)\rightarrow \Phi_{\geq -1}$
induces a bijection between the following sets:

 (1) The set of basic tilting objects in $\mathcal{C}(\Omega)$;

  (2) The set of clusters in $\Phi_{\ge
-1}$.

Moreover, if we take the orientation $\Omega$ to be the $\Omega_0$
such that $(\G, \Omega_0)$ is an alternating valued quiver, then
the bijection  $\mathcal{P}\circ\gamma _{\Omega_0}$ from
$\ind\mathcal{C}(\Omega_0)$ to the set of cluster variables of a cluster algebra of type $\G$ sends basic tilting objects in
$\mathcal{C}(\Omega_0)$ to clusters of this cluster algebra, where $\mathcal{P}$ is the 1-1 correspondence from
 $\Phi_{\geq -1}$ to the set of cluster variables
of $\mathcal{A}$ (Theorem 1.9. [10]).}
\medskip

\textbf{Proof.} The subset $A$ of $\Phi _{\geq -1}$ is a cluster
if and only if the subset $\gamma _{\Omega}^{-1}(A)$ of
ind$\mathcal{C}(\Omega)$ is a basic tilting set. Combining with
Theorem 1.9 in [10], we finish the proof.
\medskip

By this theorem, we have that the tilting graph $\Delta _{\Omega}$
is a realization of exchange graph $E(\Phi)$ in [11]. Then
Theorem 5.1 in [4] gives a quiver interpretation of Theorem
1.15 [11].
\medskip

\textbf{Corollary 3.8.} {\it For every cluster $C$ and every
element $\alpha \in C$, there is a unique cluster $C'$ such that
$C\cap C'=C-{\alpha}.$ Thus the exchange graph $E(\Phi)$ is
regular of degree $n:$ every vertex in $E(\Phi)$ is incident to
precisely $n$ edges.}
 \medskip

\begin{center}
\textbf {ACKNOWLEDGMENTS.}\end{center} The author would like to thank Professor Idun Reiten for helpful
conservation on this topic. He is grateful to the
 referee for a number of helpful comments and valuable suggestions.

\begin{center}

\end{center}

\medskip
\end{document}